\documentclass[11pt]{amsart}
\usepackage[a4paper,hmargin=3cm,vmargin=3cm]{geometry}
\usepackage{amsfonts,amssymb,amscd,amstext,amsthm}
\usepackage[pdftex]{graphicx}
\usepackage[dvips]{epsfig}
\usepackage[T1]{fontenc}
\usepackage[english]{babel}
\usepackage[utf8]{inputenc}
\usepackage[bookmarksnumbered,plainpages]{hyperref}
\usepackage{indentfirst}
\usepackage{color}
\usepackage{url}

\usepackage[shortlabels]{enumitem}
\usepackage{imakeidx}
\usepackage{titlesec}
\usepackage{mathrsfs}
\usepackage{stmaryrd}
\usepackage{textfit}
\usepackage{lipsum}
\usepackage{amsmath}
\usepackage{mathtools}

\usepackage{fancyhdr}
\usepackage{times}

\def \R{\mbox{${\mathbb R}$}}
\def \S{\mbox{${\mathbb S}$}}

\newcommand{\co}{\colon}

\newcommand{\rosStyle}[1]{
    \begin{center}
        \begin{quote}
            \textit{``#1''}
        \end{quote}
    \end{center}
}

\newenvironment{proofof}[1]{
    \proof[\textbf{Proof of Theorem #1}]
}{
    \endproof
}

\usepackage{lineno}

\titleformat{\section}
{\filcenter\bfseries\large} {\thesection{.}}{0.2cm}{}
\titleformat{\subsection}[runin]
{\bfseries} {\thesubsection{.}}{0.15cm}{}[.]
\titleformat{\subsubsection}[runin]
{\em}{\thesubsubsection{.}}{0.15cm}{}[.]

\usepackage[up,bf]{caption}

\newtheorem{theorem}{Theorem}[section]

\newtheorem{thm}{Theorem}

\newtheorem{remark}[theorem]{Remark}

\theoremstyle{definition}
\numberwithin{equation}{section}
\numberwithin{figure}{section}

\begin{document}
    
\fancyhead[LO]{A congruence theorem for compact embedded hypersurfaces in \(\S^{n+1}_+\)}
\fancyhead[RE]{Allan Freitas, Felippe Guimarães}
\fancyhead[RO,LE]{\thepage}

\thispagestyle{empty}

\begin{center}
{\bf \LARGE A congruence theorem for compact embedded hypersurfaces in \(\S^{n+1}_+\)}
\vspace*{5mm}

\hspace{0.2cm} {\Large Allan Freitas and Felippe Guimarães} 
 
\end{center}



\begin{quote}
{\small
\noindent {\bf Abstract}\hspace*{0.1cm}
    We prove a codimension reduction and congruence theorem for compact \(n\)-dimensional submanifolds of \(\S^{n+p}\) that admit a mean convex isometric embedding into \(\S^{n+1}_+\) using a Reilly type formula for space forms.\\

}

{\small
\noindent {\textbf{Mathematics Subject Classification:}}\hspace*{0.1cm}
    53A07, 53C42.
}
\\
{\small
\noindent {\textbf{Keywords:}}\hspace*{0.1cm}
    Isometric embedding, Reilly type formula, Isometric rigidity.
}

\end{quote}


\section{Introduction}
    In the realm of submanifold theory, a relevant research topic is the uniqueness of isometric immersions \( f\co M^m \to N^n \) between Riemannian manifolds, where the superscript denotes the dimension. Since any isometric immersion composed with an isometry \( \tau\co N^n \to N^n \), \( g = \tau \circ f \), remains an isometric immersion, uniqueness is considered in terms of congruences by isometries of the larger space. In this context, an isometric immersion is considered isometrically rigid if all other isometric immersions are congruent to it through an isometry of the ambient space. When the ambient space is a simply connected space form, that is, a simply connected Riemannian manifold with constant sectional curvature, some results in the literature ensure local isometric rigidity under conditions on the non-degeneracy of the second fundamental form, namely \cite{KillingRigidity,AllendoerferRigidity,doCarmoDajczerRigidity,DajczerFloritGenuine} (see \cite{DajczerTojeiroBook} for a more detailed discussion). While globally the results are more scarce, here we mention the main result in \cite{SackstederRigidity} which ensures that, given two isometric immersions \(f_1, f_2: \Sigma^{n} \rightarrow \R^{n+1}\), \(n\geq 3\), from a compact hypersurface, pointwise they have the same second fundamental form, up to a sign (cf. \cite[Theorem 13.2]{DajczerTojeiroBook} for a more general statement). Global results that guarantee the congruence of the isometric immersion of a compact hypersurface in a space form, if it has constant scalar curvature, were obtained in \cite{HarleRigidity}, and results for when it is convex, in the sense that the second fundamental form is positive definite in a suitable normal direction, can be found in \cite{doCarmoWarnerConvexity} (see \cite{deLimaSurveyConvex} and references therein for a good discussion). In \cite{RosAlexandrovScalar}, the author showed that the sphere is the only compact hypersurface with constant scalar curvature embedded in Euclidean space. Additionally, a congruence result was presented, making ingenious use of Reilly's formula presented in \cite{ReillyFormula}.
    \begin{theorem}[Theorem 2 in \cite{RosAlexandrovScalar}]\label{thm:Ros}
        Let \(\Psi\co M^n \rightarrow \R^{n+1}\) be a compact embedded hypersurface. Suppose that the mean curvature of \(\Psi\), with respect to the interior normal \(H\), is nonnegative. If $\Psi'\co M^n \rightarrow \R^{n+p}$, $p\geq 1$ is an isometric immersion with mean curvature vector $\mathcal{H}$ satisfying $\Vert\mathcal{H}\Vert \leq H$ everywhere, then $\Psi' = i\circ \Psi$ where $i$ is a totally geodesic inclusion $i\co\R^{n+1}\rightarrow\R^{n+p}$. In particular, $\Psi'$ and $\Psi$ are congruent. 
    \end{theorem}
        In the last few decades, with the aim of obtaining Alexandrov-type theorems in ambient spaces beyond Euclidean space, generalizations of Reilly's formula have emerged, particularly noteworthy are those for space forms presented in \cite{QiuXiaReillySpaceForms} and more recently for sub-static manifolds as discussed in \cite{LiXiaReillySubstatic}. Specifically by employing the formulation presented in \cite{QiuXiaReillySpaceForms}, we are able to derive a congruence theorem for hypersurfaces in a hemisphere on the sphere, which is analogous to Theorem \ref{thm:Ros}.
    \begin{thm}\label{thm:MainA}
        Let \(\Psi\co M^n \rightarrow \S^{n+1}_+\) be a compact hypersurface embedded into $\S^{n+1}_+$. Suppose that the mean curvature of \(\Psi\), with respect to the interior normal \(H\), is nonnegative. If $\Psi'\co M^n \rightarrow \S^{n+p}$, $p\geq 1$ is an isometric immersion with mean curvature vector $\mathcal{H}$ satisfying $\Vert\mathcal{H}\Vert \leq H$ everywhere, then $\Psi' = i\circ \Psi$ where $i$ is a totally geodesic inclusion $i\co\S^{n+1}\rightarrow\S^{n+p}$. In particular, $\Psi'$ and $\Psi$ are congruent. 
    \end{thm}

    \begin{remark}\label{rmk:OverHyp}
        Observe that in Theorem \ref{thm:MainA} we do not assume embedding and that the image for $\Psi'$ is contained within a hemisphere. It is worth mentioning that the hypersurface must have some points with strictly positive mean curvature, as we know that a minimal hypersurface cannot be contained within a hemisphere.
    \end{remark} 
    \noindent As discussed in \cite{RosAlexandrovScalar}, this result can be seen as an extension of the classic Schur's theorem to hypersurfaces in a hemisphere of the sphere. The case \(n=1\) was addressed in a more general setting in \cite{NiSchur} using different techniques. Furthermore, this type of result is naturally connected to the question:
    \rosStyle{Are two isometric immersions \(f_1, f_2: \Sigma^{n} \rightarrow M^{n+1}\) with same mean curvature necessarily congruent?}
    
   \noindent A good discussion about this question can be found in \cite{LiMiaoIsoRigiditySameH} and references therein, as well as a related result for hypersurfaces in a spatial Schwarzschild or AdS-Schwarzschild manifold. Integral formulas have been important for achieving results on intrinsic rigidity, specifically regarding the extension of geometric properties from a manifold which is the boundary of a region to the region itself. This investigation is distinct from the main question of this article, yet the techniques show similarities, as showed in \cite{HangWangNonnegativeRicci}, where the approach introduced in \cite{RosAlexandrovScalar} is used to prove the following result.
    \begin{theorem}(Proposition 2 in \cite{HangWangNonnegativeRicci})\label{thm:HangWang}
        Let \(M^n\) be a smooth compact connected Riemannian manifold with connected boundary \(\Sigma = \partial M\) and \(Rc \geq 0\). If \(\phi: \Sigma \rightarrow \R^n\) is an isometric immersion with \(\Vert\mathbf{H}_\phi\Vert \leq H\) on \(\Sigma\), where \(\mathbf{H}_\phi\) is the mean curvature vector of the immersion \(\phi\), then \(M\) is flat. If \(\phi\) is also an embedding, \(M\) is isometric to a domain in \(\R^n\).
    \end{theorem}
     \noindent This is a generalization of Theorem \ref{thm:Ros} by requiring that the ambient space has non-negative Ricci curvature, which leads to the conclusion that the region must be flat. Further research related to this result was conducted in \cite{HangWangNonnegativeRicci,LaiNote,MiaoWangConjecture}, placing a lower bound on Ricci curvature compared to space forms. However, unlike Theorem \ref{thm:HangWang}, which required the boundary to be, in a sense, mean convex, the results in \cite{MiaoWangConjecture} were achieved under a stronger hypothesis of convexity, similar to those in \cite{doCarmoWarnerConvexity}. The authors also propose the following question:
     \rosStyle{Let \((\Omega,g)\) be an \(n\)-dimensional, compact Riemannian manifold with boundary \(\Sigma\). Let \(D \subset S^n_+\) be a bounded domain with smooth boundary \(\partial D\), where \(S^n_+\) is the standard n-dimensional hemisphere. Suppose the following
     \begin{itemize}
        \item Ric \(\geq (n-1)g\).
        \item There exists an isometry \(X: \Sigma \rightarrow \partial D\).
        \item \(H \geq H_{S} \circ X\), where \(H, H_{S}\) are the mean curvature of \(\Sigma, \partial D\) in \((\Omega,g), S^n_+\) respectively.
     \end{itemize}
     Is \((\Omega,g)\) isometric to \(D\) in \(S^n_+\)?}
    A positive answer to this question would imply that our theorem could be improved in a similar way to Theorem \ref{thm:HangWang}.

\section{The proof}\label{sec:Proofs}
In the proof, we utilize a generalization of Reilly's formula by Qiu and Xia (cf. \cite[Theorem 1.1]{QiuXiaReillySpaceForms}), which particularly yields an intriguing identity in the spherical case:
    \begin{equation}\label{eq:ReillySpaceforms}
        \begin{aligned}
            \int_{\Omega} \cos r(x)((\bar{\Delta} f + (n+1)f)^2 - |\bar{\nabla}^2 f + fg|^2) d\Omega
            &= \int_{M} \cos r(x)(2u \Delta z + nHu^2 + h(\nabla z, \nabla z)) \\
            &\quad + 2n uz) dA \\
            &\quad + \int_{M} \overline{\nabla}_{\nu} [\cos r(x)](|\nabla z|^2 - n z^2) dA.
        \end{aligned}
    \end{equation}
Here $\Omega \subset \S^{n+1}_{+}$ is a domain with boundary $M$, $r(x) = d_{g}(x, p)$ where $p \in \S^{n+1}_+$ is the corresponding pole of the hemisphere, and $f \co \Omega \rightarrow \R$ is a smooth function. Furthermore, $\bar{\nabla}$ and $\bar{\Delta}$ indicate the gradient and the Laplacian operators in $\S^{n+1}_{+}$,  $z=f\vert_{M}$ and $u=\bar{\nabla}_{\nu}f$, where $\nu$ is the unit outward normal of $M$. Furthermore, $\Delta$ and $\nabla$ denote the correspondent intrinsic operators in $M$, $h(X, Y)= g(\bar{\nabla}_{X}\nu, Y)$ and $nH=tr_{g}h$ are the second fundamental form and the mean curvature (with respect to $\nu$) of $M$, respectively.

\begin{proofof}{\ref{thm:MainA}}
   Let \(\Tilde{\Psi} = \mathcal{I}\circ \Psi'\) where \(\mathcal{I}: \S^{n+p} \rightarrow \R^{n+p+1}\) is the umbilical inclusion. Choose a global parallel frame \(\{ e_1,e_2,\ldots,e_{n+p+1}\}\) of \(\R^{n+p+1}\) and define \(\Tilde{\Psi}_i: M^{n} \rightarrow \R\) given by \(\Tilde{g}(\Tilde{\Psi},e_i)\) where \(\Tilde{g}\) is the standard metric of \(\R^{n+p+1}\). Also, consider $\Omega\subset\S^{n+1}_{+}$ the bounded domain such that $\partial\Omega=M$.
    \noindent Before using \eqref{eq:ReillySpaceforms} for a suitable function $f$, let's simplify some terms. Let \(\{v_1,...,v_n\}\) be a local orthonormal frame in the tangent bundle of \(M^n\), and using that $\Psi$ is an isometric immersion, we have that
    \[g(\nabla \Tilde{\Psi}_i, v_s) = v_s \Tilde{g}(\Tilde{\Psi}, e_i) = \Tilde{g}(\Tilde{\Psi}_* v_s, e_i),\] which also implies that \(\Tilde{\Psi}_*\nabla \Tilde{\Psi}_i = e_i - \sum^{p+1}_{k=1}\Tilde{g}(e_i, \xi_k)\xi_k\), where \(\{ \xi_1, \xi_2, \ldots, \xi_{p+1}\}\) is a local orthonormal frame in the normal bundle of \(M^n\) with respect to the isometric immersion \(\Tilde{\Psi}\). In particular 
    \[\sum^{n+p+1}_{i=1} |\nabla \Tilde{\Psi}_i|^2 = \sum^{n+p+1}_{i=1}\left(1 - \sum^{p+1}_{k=1}\Tilde{g}(e_i, \xi_k)^2\right) = n,\] and 
    
    \begin{equation}\label{eq1}
    \begin{split}    
    \sum_{i=1}^{n+p+1} h(\nabla \tilde{\Psi}_i, \nabla\tilde{\Psi}_i) & = \sum_{i=1}^{n+p+1}\sum_{s,k=1}^{n} \tilde{g}(\tilde{\Psi}_*v_s,e_i)\tilde{g}(\tilde{\Psi}_*v_k,e_i) h(v_s,v_k)\\ & = \sum_{s,k=1}^{n} g_{sk} h(v_s,v_k) = nH.
    \end{split}
    \end{equation}
    As \(\Tilde{\Psi}(M^n)\subset \S^{n+p}\subset \R^{n+p+1}\) it follows that $\Vert \tilde{\Psi}(x) \Vert^2=1$ for all \(x\in M^n\). This implies 
    
    \begin{equation}\label{eq2}
     \sum^{n+p+1}_{i=1} \left[|\nabla \Tilde{\Psi}_i|^2 -n(\tilde{\Psi}_i)^2\right] = 0.   
    \end{equation}
    
    Now, for $1\leq i\leq n+p+1$, let $f_i\co \overline{\Omega}\rightarrow \R$ be the solution of the problem
    \[
      \begin{cases}
                   \Delta f_i + (n+1)f_i = 0 \ &\quad \text{in} \quad\Omega\\
                   f_i = \tilde{\Psi}_i&\quad \text{on}\quad M 
                \end{cases}
    \]

    \noindent By considering \(u_i = \frac{\partial f_i}{\partial \nu}\) and $U=(u_{1},\ldots,u_{n+p+1})$, using \eqref{eq1} and $\Delta\Tilde{\Psi}=n\mathcal{H}-n\Tilde{\Psi}$, we have:
    
    \begin{equation}\label{eq3}
        \begin{aligned}
            \sum^{n+p+1}_{i=1} & \left[2u_i\Delta \tilde{\Psi}_i +nHu_i^2+h(\nabla \tilde{\Psi}_i, \nabla \tilde{\Psi}_i) +2n u_i\tilde{\Psi}_i )\right]\\ 
            &= \left[ 2\tilde{g}(U, \Delta \tilde{\Psi}) + nH\Vert U \Vert^2 +nH +2n\tilde{g}(U, \tilde{\Psi}) \right] \\
            &= \left[ 2\tilde{g}(U, n\mathcal{H}-n\tilde{\Psi}) + nH\Vert U \Vert^2 +nH +2n\tilde{g}(U, \tilde{\Psi}) \right] \\
            &= n\left[H\Vert U \Vert^2 + 2\Tilde{g}(U,\mathcal{H}) +H \right].
        \end{aligned}
    \end{equation}
    
     \noindent Applying the Reilly's identity \eqref{eq:ReillySpaceforms} for each $f_{i}$, summing up from $1$ to $n+p+1$ and using \eqref{eq2} and \eqref{eq3} we get the following expression \[-\int_{\Omega} \cos r(x)\sum^{n+p+1}_{i=1}(|\nabla^2 f_i - \frac{1}{n+1}(\Delta f_i)g |^2) d\Omega
            = \int_{M} n\cos r(x) \left[H\Vert U \Vert^2 + 2\Tilde{g}(U,\mathcal{H}) +H \right].\]
    Moreover, using that \(2\Tilde{g}(U,\mathcal{H}) \geq - 2 \Vert U \Vert  \Vert\mathcal{H}\Vert \geq - 2 \Vert U \Vert  H\) it follows \[0\geq-\int_{\Omega} \cos r(x)\sum^{n+p+1}_{i=1}(|\nabla^2 f_i - \frac{1}{n+1}(\Delta f_i)g |^2) d\Omega \geq \int_{M} H\cos r(x) \left[\Vert U \Vert  - 1 \right]^2.\] Hence, \(\Vert U \Vert =1\) in every point where $H\neq 0$, and \(|\nabla^2 f - \frac{1}{n+1}(\Delta f_i)g | = 0\) for all \(1\leq i \leq n+p+1\). This tells us that
    \begin{equation}\label{eq:Conclusion}    
    \begin{split}
        \Vert U \Vert =1\ \text{at } x \in M &\Rightarrow \Vert\Tilde{\Psi}_*X\Vert = |X| \ \text{for all }X\in TM\oplus\text{span}\{\nu\};\\
        |\nabla^2 f_i - \frac{1}{n+1}(\Delta f_i)g | = 0 \ \text{at } x \in \Omega &\Rightarrow \nabla^2 f_i = - f_ig.
    \end{split}
    \end{equation}
    The rest of the proof can be addressed in a more general manner, but here we will present it for the specific setting of this work, and it can be seen as an argument similar to Obata's rigidity. Let \(f = (f_1, f_2, \ldots, f_{n+p+1})\co \Omega \to \R^{n+p+1}\) and define \(A = \{ x \in \Omega \mid \Vert f(x)\Vert^2 = 1 \text{ and } \Vert f_*(x)(w)\Vert = |w| \text{ for all } w \in T_x \Omega \}\). Due to the existence of a point $x \in M$ with \(H(x) \neq 0\) (which exists in view of Remark \ref{rmk:OverHyp}) and from \eqref{eq:Conclusion}, we deduce that \(A \neq \emptyset\). Given the smoothness of \(f\), \(A\) is closed. In order to prove that $A=\Omega$ remaining to show that \(A\) is open. First, the solutions of the equations \(\nabla^2 f_i = - f_ig\) for all \(1 \leq i \leq n+p+1\) can be extended to an open neighborhood \(\Omega^a\) of \(\Omega\), thus avoiding dealing with boundary points. In order to solve the equation in \eqref{eq:Conclusion} explicitly, let \(y \in A\) and \(U \subset \Omega^a\) be a totally normal neighborhood of \(y\). It follows from \(\nabla^2 f_i = - f_ig\) for all \(1 \leq i \leq n+p+1\) that
    \begin{equation}\label{eq:formaGeralf}
        f(x) = \cos(r_y(x)) f(y) + \sin(r_y(x)) f(y)_*(\nabla r_y(x)),
    \end{equation}
    where \(r_y(x) = d(x,y)\) for \(x \in \Omega\). From \(y \in A\), we know that \(\Vert f(y) \Vert^2=1\) holds for every \(y \in U\). Now considering \(x \in U \backslash \{y\}\), \(w \in T_x\Omega^a\) with \(\Vert w\Vert=1\), and \(\gamma: [0,d(x,y)]\rightarrow U\) as the geodesic connecting \(\gamma(0) = x\) to \(\gamma(d(x,y)) = y\), let \(\alpha: [0,T) \rightarrow U\) be a unit geodesic such that \(\alpha(0)=x\) and \(\alpha'(0) = w\). First, assuming that \(\alpha = \gamma\), we observe that
        \[f_*(x)w=\frac{d}{ds}f\circ\alpha(s)\biggr\rvert_{s = 0} = \sin(r_y(x))f(y) - \cos(r_y(x))f_*(y)w \]
    as in this case \(\nabla r_y(\alpha(s)) = -\alpha'(d(x,y))\), and since \(y \in A\) it implies that \(\Vert f_*(x)(w)\Vert^2 = 1 = |w|^2\). We can now assume that \(g(w, \gamma'(0)) = 0\). By using the exponential map of \(\S^{n+1}_+ \subset \R^{n+2}\) at a point \(y\), we have that
        \[\alpha(s) = \cos r_y(\alpha(s))y + \sin{r_y(\alpha(s))} \nabla r_y(\alpha(s)).\]
    Hence
        \[\frac{d}{ds} \nabla r_y(\alpha(s))\biggr\rvert_{s = 0}= \frac{w}{\sin{r_y(x)}}.\]
    This leads to the conclusion that
        \[f_*(x)w=\frac{d}{ds}f\circ\alpha(s)\biggr\rvert_{s = 0} = \sin(r_y(x))f_*(y) \frac{d}{ds} \nabla r_y(\alpha(s))\biggr\rvert_{s = 0} = f_*(y)w,\]
    and given that \(y \in A\) it follows that $U \subset A$, concluding that \(A = \Omega\). Here, by using the inclusion of the sphere in Euclidean space, we are clearly assuming some identifications, since \(f_*(x)w = \sin(r_y(x))f_*(y)\left[{exp_y}_*(\frac{d}{ds} \nabla r_y(\alpha(s))\rvert_{s = 0})\right]^{-1}w\). Consequently, \(f\) is an isometric immersion, and as \(\Omega \subset \S^{n+1}_+\), \(f\) can be described in the entire hemisphere using an exponential chart at a point \(x \in \Omega\) as indicated in \eqref{eq:formaGeralf}. This expression ensures that \(f(\Omega) \subset \S^{n+1} \subset \S^{n+p}\), where the totally geodesic submanifold \(\S^{n+1}\) is determined as the sphere centered at the origin with its tangent space given by \(f_*(x)(T_x\Omega)\). Given the boundary condition that \(f \circ \Psi = \mathcal{I} \circ \Psi'\), this concludes the proof of the theorem.
\end{proofof}

    In the proof presented, with regard to the hypothesis about the mean curvature of the isometric immersions, we are actually using the following inequality
    \[
    \int_{\Omega} \cos(r(x))H d\Omega \geq \int_{\Omega} \cos(r(x))\Vert \mathcal{H}\Vert d\Omega.
    \]
    Thus, we can rewrite our result as follows.
    \begin{thm}\label{thm:MainB}
        Let \(\Psi: M^n \rightarrow \S^{n+1}_+\) be a compact hypersurface embedded into \(\S^{n+1}_+\). Suppose that the mean curvature \(H\) of \(\Psi\), with respect to the interior domain \(\Omega\), is nonnegative. Let \(\Psi': M^n \rightarrow \S^{n+p}\), \(p\geq 1\) be an isometric immersion then we have that the following inequality is satisfied \[\int_{\Omega} \cos(r(x))H d\Omega \leq \int_{\Omega} \cos(r(x))\Vert \mathcal{H}\Vert d\Omega,\] where \(\mathcal{H}\) is the mean curvature vector of \(\Psi'\) and $r(x) = d_{g}(x, p)$ where $p \in \S^{n+1}_+$ is the corresponding pole of the hemisphere. Moreover, if equality holds then \(\Psi' = i\circ \Psi\) where \(i\) is a totally geodesic inclusion \(i: \S^{n+1}\rightarrow \S^{n+p}\). In particular, \(\Psi'\) and \(\Psi\) are congruent.
    \end{thm}

\section*{Acknowledgements}
The authors would like to thank Theodoros Vlachos and Ronaldo Lima for their valuable suggestions concerning the manuscript.

Allan Freitas has been partially supported by CNPq/Brazil (grant 316080/2021-7) and Programa Primeiros Projetos - FAPESQ/PB and MCTIC/CNPq (grant 2021/3175). Felippe Guimar\~aes is supported by FAPESQ/PB and partially by CNPq/Brazil (grant 409513/2023-7).





\bibliographystyle{abbrv}
\bibliography{bibliography}

\begin{thebibliography}{10}

\bibitem{AllendoerferRigidity}
C.~B. Allendoerfer.
\newblock Rigidity for spaces of class greater than one.
\newblock {\em Am. J. Math.}, 61:633--644, 1939.

\bibitem{DajczerFloritGenuine}
M.~Dajczer and L.~A. Florit.
\newblock Genuine deformations of submanifolds.
\newblock {\em Commun. Anal. Geom.}, 12(5):1105--1129, 2004.

\bibitem{DajczerTojeiroBook}
M.~Dajczer and R.~Tojeiro.
\newblock {\em Submanifold theory. Beyond an introduction}.
\newblock Universitext. Springer, New York, 2019.

\bibitem{deLimaSurveyConvex}
R.~de~Lima.
\newblock A survey on convex hypersurfaces of riemannian manifolds.
\newblock {\em Matemática Contempor\^anea}, 50(8), 2022.

\bibitem{doCarmoDajczerRigidity}
M.~do~Carmo and M.~Dajczer.
\newblock A rigidity theorem for higher codimensions.
\newblock {\em Math. Ann.}, 274(4):577--583, 1986.

\bibitem{doCarmoWarnerConvexity}
M.~P. do~Carmo and F.~W. Warner.
\newblock {Rigidity and convexity of hypersurfaces in spheres}.
\newblock {\em Journal of Differential Geometry}, 4(2):133 -- 144, 1970.

\bibitem{HangWangNonnegativeRicci}
F.~Hang and X.~Wang.
\newblock Vanishing sectional curvature on the boundary and a conjecture of
  {Schroeder} and {Strake}.
\newblock {\em Pac. J. Math.}, 232(2):283--287, 2007.

\bibitem{HarleRigidity}
C.~E. Harle.
\newblock Rigidity of hypersurfaces of constant scalar curvature.
\newblock {\em J. Differ. Geom.}, 5:85--111, 1971.

\bibitem{KillingRigidity}
W.~Killing.
\newblock Analytical treatment of non-euclidean spaces.
\newblock Leipzig. {Teubner}. {XII} u. 264 {Seiten} (1885)., 1885.

\bibitem{LaiNote}
M.~Lai.
\newblock A note on {Hang}-{Wang}'s hemisphere rigidity theorem.
\newblock {\em Math. Z.}, 296(3-4):901--909, 2020.

\bibitem{LiMiaoIsoRigiditySameH}
C.~Li, P.~Miao, and Z.~Wang.
\newblock Uniqueness of isometric immersions with the same mean curvature.
\newblock {\em J. Funct. Anal.}, 276(9):2831--2855, 2019.

\bibitem{LiXiaReillySubstatic}
J.~Li and C.~Xia.
\newblock An integral formula and its applications on sub-static manifolds.
\newblock {\em J. Differ. Geom.}, 113(3):493--518, 2019.

\bibitem{MiaoWangConjecture}
P.~Miao and X.~Wang.
\newblock Boundary effect of {Ricci} curvature.
\newblock {\em J. Differ. Geom.}, 103(1):59--82, 2016.

\bibitem{NiSchur}
L.~Ni.
\newblock A {Schur}’s theorem via a monotonicity and the expansion module.
\newblock {\em Journal für die reine und angewandte Mathematik (Crelles
  Journal)}, 2023(805):101--114, 2023.

\bibitem{QiuXiaReillySpaceForms}
G.~Qiu and C.~Xia.
\newblock A generalization of {Reilly}'s formula and its applications to a new
  {Heintze}-{Karcher} type inequality.
\newblock {\em Int. Math. Res. Not.}, 2015(17):7608--7619, 2015.

\bibitem{ReillyFormula}
R.~C. Reilly.
\newblock Applications of the hessian operator in a riemannian manifold.
\newblock {\em Indiana University Mathematics Journal}, 26(3):459--472, 1977.

\bibitem{RosAlexandrovScalar}
A.~Ros.
\newblock Compact hypersurfaces with constant scalar curvature and a congruence
  theorem.
\newblock {\em J. Differ. Geom.}, 27(2):215--220, 1988.

\bibitem{SackstederRigidity}
R.~Sacksteder.
\newblock The rigidity of hypersurfaces.
\newblock {\em J. Math. Mech.}, 11:929--939, 1962.

\end{thebibliography}

\vskip 0.2cm

\noindent Allan Freitas

\noindent Departamento de Matematica, \\ Universidade Federal da Paraíba, João Pessoa (Brazil).

\noindent  e-mail: {\tt allan@mat.ufpb.br}

\vskip 0.2cm

\noindent Felippe Guimarães

\noindent Departamento de Matematica, \\ Universidade Federal da Paraíba, João Pessoa (Brazil).

\noindent  e-mail: {\tt fsg@academico.ufpb.br}

\end{document}